\documentclass{article}
\usepackage{amssymb}
\usepackage{amsfonts}
\usepackage{amsmath}

\input{tcilatex}
\begin{document}

\title{On the Kurzweil-Henstock Integral in Probability}
\author{Sorin G. Gal \\
Department of Mathematics and Computer Science, \\
University of Oradea, \\
Universitatii 1, 410087, Oradea, Romania\\
E-mail: \textit{galso@uoradea.ro}}
\date{}
\maketitle

\begin{abstract}
By using the method in [5], the aim of the present note is to generalize the Riemann integral in probability introduced in [7],
to Kurzweil-Henstock integral in probability. Properties of the new integral are proved.
\end{abstract}



\section{Introduction}

Let $(E, B, P)$ be a field of probability, where $E$ is a nonempty set, $B$ a field of parts on $E$ and $P$ a composite probability on $B$.
Let us denote by $L(E, B, P)$ the set of all real random variables defined on $E$ and a.e. finite.

It is well-known the following concept :

{\bf Definition 1.1.} (see e.g. [7], p. 50) We say that the random function $f:[a, b]\to L(E,B,P)$ (where $a, b\in \mathbb{R}, a <b$) is Riemann integrable in probability on $[a, b]$, if there exists a random variable $I=I(\omega)\in L(E,B,P)$ satisfying :
for all $\varepsilon > 0, \eta >0$, there exists $\delta=\delta(\varepsilon)>0$, such that for all divisions $d: a=x_{0}<x_{1}< ... <x_{n}=b$ with
the norm $\nu(d)<\delta$ and all $\xi_{i}\in [x_{i}, x_{i+1}]$, $i\in \{0, ..., n-1\}$, we have
$$P(\{\omega\in E ; |S(f;d, \xi_{i})(\omega)-I(\omega)|\ge \varepsilon\})<\eta,$$
where $\nu(d)=\max\{x_{i+1}-x_{i} ; i=0, 1, ..., n-1\}$ and
$$S(f;d, \xi_{i})(\omega)=\sum_{i=0}^{n-1}f(\xi_{i}, \omega)(x_{i+1}-x_{i}).$$
In this case, $I(\omega)$ is called the Riemann integral in probability of $f$ on $[a, b]$ and it is denoted by $I(\omega)=(P)\int_{a}^{b}f(t, \omega)dt$.

{\bf Remark.} As it was proved in [7, p. 50], if $I_{1}(\omega), I_{2}(\omega)$ are Riemann integrals in probablity of $f$ on $[a, b]$,
then $P(\{\omega\in E ; I_{1}(\omega)\not= I_{2}(\omega)\})=0$.

Using the method in [5], in Section 2 we introduce the so-called Kurzweil-Henstock (Riemann generalized) integral in probability. Section 3 contains basic properties of this generalized integral.

\section{The Kurzweil-Henstock Integral in Probability}

Firstly, we recall some concepts in [5] we need for our purpose.

A tagged division of $[a, b]$ is of the form
$$d_{S}: a=x_{0}\le \xi_{0}\le x_{1}< ... <x_{i}\le \xi_{i}\le x_{i+1}< ... < x_{n-1}\le \xi_{n-1}\le x_{n}=b.$$
A gauge on $[a, b]$ is an open interval valued function $\gamma$ defined on $[a, b]$, such that $t\in \gamma(t)$, $t\in [a, b]$.
A tagged division $d_{S}$ of $[a, b]$ is called $\gamma$-sharp if $[x_{i}, x_{i+1}]\subset \gamma(\xi_{i})$, for all $i\in \{0, ..., n-1\}$.

Now, we are in position to introduce the following.

{\bf Definition 2.1.} Let $f:[a, b]\to L(E, B, P)$. A random variable $I(\omega)\in L(E,B,P)$ is called a Kurzweil-Henstock (shortly (KH)) integral in probability of $f$ on $[a, b]$, if : for all $\varepsilon > 0, \eta >0$, there exists $\gamma_{\varepsilon, \eta}$-gauge on $[a, b]$, such that for any tagged division $d_{S}$ which is $\gamma_{\varepsilon, \eta}$-sharp, we have
$$P(\{\omega\in E ; |S(f;d_{S}, \xi_{i})(\omega)-I(\omega)|\ge \varepsilon\})<\eta,$$
where $S(f;d_{S}, \xi_{i})(\omega)=\sum_{i=0}^{n-1}f(\xi_{i}, \omega)(x_{i+1}-x_{i}).$

In this case, we write $I(\omega)=(KH)\int_{a}^{b}f(t, \omega)dt$.

{\bf Theorem 2.2.} {\it If $I_{1}(\omega)$, $I_{2}(\omega)$ are (KH) integrals in probability of $f$ on $[a, b]$, then $P(\{\omega\in E ; I_{1}(\omega)\not= I_{2}(\omega)\})=0$.}

{\bf Proof.} We will prove that $P(\{\omega\in E ; |I_{1}(\omega)-I_{2}(\omega)|>0\})=0$.  Let $\varepsilon, \eta >0$. There exists the gauges
$\gamma^{(1)}_{\varepsilon, \eta}$, $\gamma^{(2)}_{\varepsilon, \eta}$ on $[a, b]$, such that for any tagged divisions of $[a, b]$, $d^{(1)}_{S}$,
$d^{(2)}_{S}$ which are $\gamma^{(1)}_{\varepsilon, \eta}$-sharp and $\gamma^{(2)}_{\varepsilon, \eta}$-sharp, respectively, we have
$$P(\{\omega\in E ; |S(f;d_{S}^{(1)}, \xi_{i})(\omega)-I_{1}(\omega)|\ge \varepsilon/2\})<\eta/2,$$
$$P(\{\omega\in E ; |S(f;d_{S}^{(2)}, \xi_{i})(\omega)-I_{2}(\omega)|\ge \varepsilon/2\})<\eta/2.$$
Let us define a new gauge on $[a, b]$ by $\gamma(t)=\gamma^{(1)}_{\varepsilon, \eta}(t)\bigcap \gamma^{(2)}_{\varepsilon, \eta}(t)$,
$t\in [a, b]$. By [5, Section 1.8], there exists a $\gamma$-sharp tagged division $d_{S}$ of $[a, b]$.

Since $\gamma(t)\subset \gamma^{(1)}_{\varepsilon, \eta}(t)$, $\gamma(t)\subset \gamma^{(2)}_{\varepsilon, \eta}(t)$, $t\in [a, b]$, obviously
that $d_{S}$ is $\gamma^{(1)}_{\varepsilon, \eta}$-sharp and $\gamma^{(2)}_{\varepsilon, \eta}$-sharp too.

We have
$$|I_{1}(\omega)-I_{2}(\omega)|\le |I_{1}(\omega)-S(f;d_{S}, \xi_{i})(\omega)|+|S(f; d_{S}, \xi_{i})(\omega)-I_{2}(\omega)|,$$
which immediately implies
$$\{\omega\in E ; |I_{1}(\omega)-I_{2}(\omega)|\ge \varepsilon\}\subset \{\omega\in E; |I_{1}(\omega)-S(f;d_{S}, \xi_{i})(\omega)|\ge \varepsilon/2\}$$
$$\bigcup \{\omega\in E; |S(f;d_{S}, \xi_{i})(\omega)-I_{2}(\omega)|\ge \varepsilon/2\}$$
and
$$P(\{\omega\in E ; |I_{1}(\omega)-I_{2}(\omega)|\ge \varepsilon\})\le P(\{\omega\in E; |I_{1}(\omega)-S(f;d_{S}, \xi_{i})(\omega)|\ge \varepsilon/2\})$$
$$+ P(\{\omega\in E; |S(f;d_{S}, \xi_{i})(\omega)-I_{2}(\omega)|\ge \varepsilon/2\})<\eta/2 + \eta/2=\eta.$$
Now, considering $\varepsilon >0$ fixed and passing to limit with $\eta\to 0$, we get $P(\{\omega\in E ; |I_{1}(\omega)-I_{2}(\omega)|\ge \varepsilon\})=0$.

For $\varepsilon=\frac{1}{n}$, let us denote $A_{n}=\{\omega\in E ; |I_{1}-I_{2}|\ge 1/n\}$. Obviously $A_{n}\subset A_{n+1}$ and
$\bigcup_{n=1}^{\infty}A_{n}=\{\omega\in E ; |I_{1}(\omega)-I_{2}(\omega)|>0\}$. Then,
$$P(\{\omega\in E ; |I_{1}(\omega)-I_{2}(\omega)|>0\})=\lim_{n\to \infty}P(A_{n})=0,$$
which proves the theorem.  $\hfill \square$

As in the case of usual real functions, another definition for the (KH) integral can be the following.

{\bf Definition 2.3.} Let $f:[a, b]\to L(E,B,P)$. We say that $f$ is Kurzweil-Henstock integrable in probability on $[a, b]$, if there exists
$I\in L(E,B,P)$ with the property : for all $\varepsilon > 0$ and $\eta > 0$, there exists $\delta_{\varepsilon, \eta}:[a, b]\to \mathbb{R}_{+}$,
such that for any division $d_{S}: a=x_{0}<x_{1}< ... < x_{n}=b$ and any $\xi_{i}\in [x_{i}, x_{i+1}]$ with $x_{i+1}-x_{i} < \delta_{\varepsilon, \eta}(\xi_{i})$, $i=0, ..., n-1$, we have
$$P(\{\omega\in E ; |S(f; d_{S}, \xi_{i})(\omega)-I(\omega)|\ge \varepsilon\})<\eta.$$
{\bf Remarks.} 1) The Definitions 2.1 and 2.3 are equivalent. Indeed, this easily follows from the fact that any function $\delta_{\varepsilon, \eta}:[a, b]\to \mathbb{R}_{+}$, generates the gauge $\gamma_{\varepsilon, \eta}(t)=(t-\delta_{\varepsilon, \eta}(t)/2, t+\delta_{\varepsilon, \eta}(t)/2)$, $t\in [a, b]$ and conversely, any gauge $\gamma_{\varepsilon, \eta}$ on $[a, b]$ (which obviously can be written in the form
$\gamma_{\varepsilon, \eta}(t)=(t-\alpha(t), t+\beta(t))$, $\alpha(t), \beta(t)>0$, $t\in [a, b]$) generates the function $\delta_{\varepsilon, \eta}(t)=\alpha(t)+\beta(t)$, $t\in [a, b]$, such that the (KH)-integrability which uses
the function $\delta_{\varepsilon, \eta}$ is equivalent with the (KH)-integrability which uses the gauge $\gamma_{\varepsilon, \eta}$.

2) If $\delta_{\varepsilon, \eta}$ is a constant function, Definition 2.3 reduces to Definition 1.1.

\section{Properties of the (KH)-Integral in Probability}

In this section, we will prove some properties of the (KH)-integral in probability. Firstly, we need the following.

{\bf Definition 3.1.} (see e.g. [3, p. 82], [4]). We say that $\varphi : [a, b]\to \mathbb{R}$ is Kurzweil-Henstock integrable on $[a, b]$,
if there exists $I\in \mathbb{R}$, such that for all $\varepsilon >0$, there exists $\delta_{\varepsilon}:[a, b]\to \mathbb{R}$, such that
for any division $d: a=x_{0} < x_{1} <... < x_{n}=b$ and any $\xi_{i}\in [x_{i}, x_{i+1}]$ with $x_{i+1}-x_{i}<\delta(\xi_{i})$, we have
$|I-\sum_{i=0}^{n-1}\varphi(\xi_{i})(x_{i+1}-x_{i})|<\varepsilon$. We write $I=(KH)\int_{a}^{b}\varphi(t)dt$.

The following result holds.

{\bf Theorem 3.2.} {\it If $f:[a, b]\to L(E,B,P)$ is of the form $f(t, \omega)=\sum_{k=1}^{p}C_{k}(\omega)\cdot \varphi_{k}(t)$, where
$C_{k}\in L(E,B,P)$ and $\varphi_{k}$ are Kurzweil-Henstock integrable on $[a, b]$, $k=1, ..., p$, then $f$ is Kurzweil-Henstock integrable
in probability on $[a, b]$ and we have
$$(KH)\int_{a}^{b}f(t, \omega)dt=\sum_{k=1}^{p}C_{k}(\omega)\cdot (KH)\int_{a}^{b}\varphi_{k}(t)dt.$$}

{\bf Proof.} Obviously that it is sufficient to consider only the case when $f(t, \omega)=C(\omega)\cdot \varphi(t)$, $\omega\in E$,
$t\in [a, b]$.

If $d: a=x_{0}< ... < x_{n}=b$, $\xi_{i}\in [x_{i}, x_{i+1}]$, $i=0, ..., n-1$, then it is easy to see that $S(f;d, \xi_{i})(\omega)=C(\omega)\cdot
\sum_{i=0}^{n-1}\varphi(\xi_{i})\cdot (x_{i+1}-x_{i})$. Let us denote $I=(KH)\int_{a}^{b}\varphi(t)dt$ and $A_{m}=\{\omega\in E ; |C(\omega)|\ge m\}$. Obviously, $A_{m+1}\subset A_{m}$, $m\in \mathbb{N}$. Denoting $A=\bigcap_{m=1}^{\infty}A_{m}$, since $C\in L(E, B, P)$ we get  $P(A)=0$ and
$\lim_{m\to \infty}P(A_{m})=P(A)=0$. Consequently, if $\eta >0$, there exists $N(\eta)\in \mathbb{N}$, such tat
$$P(\{\omega\in E; |C(\omega)|\ge m\})<\eta, m\in \mathbb{N}, m\ge N(\eta).$$
For fixed $m\ge N(\eta)$, let us consider $\varepsilon >0$, such that $1/\varepsilon \ge m$. Now, for $\varepsilon^{2}>0$, since $\varphi$ is Kurzweil-Henstock integrable on $[a, b]$, by Definition 3.1, there exists $\delta_{\varepsilon^{2}}:[a, b]\to \mathbb{R}$, such that for any division
$d: a=x_{0}< ... < x_{n}=b$ and any $\xi_{i}\in [x_{i}, x_{i+1}]$ with $x_{i+1}-x_{i}<\delta_{\varepsilon^{2}}(\xi_{i})$, $i=0, ..., n-1$, we have
$|I-\sum_{i=0}^{n-1}\varphi(\xi_{i})(x_{i+1}-x_{i})| < \varepsilon^{2}$.

We have
$$\{\omega\in E; |S(f;d, \xi_{i})(\omega)-C(\omega)\cdot I|\ge \varepsilon\}$$
$$=\{\omega\in E; |C(\omega)|\cdot |I - \sum_{i=0}^{n-1}\varphi(\xi_{i})(x_{+1}-x_{i})|\ge \varepsilon\}$$
$$\subset \{\omega\in E; |C(\omega)|\ge 1/\varepsilon\}\subset \{\omega\in E ; |C(\omega)|\ge m\},$$
i.e. $P(\{\omega\in E; |S(f;d, \xi_{i})(\omega)-C(\omega)\cdot I|\ge \varepsilon\})<\eta$, for any division $d:a=x_{0}< ... <x_{n}=b$
and any $\xi_{i}\in [x_{i}, x_{i+1}]$ with $x_{i+1}-x_{i}<\delta_{\varepsilon^{2}}(\xi_{i})$ (in fact, $\varepsilon$ depends on $m$, which depends on $\eta$, therefore $\delta_{\varepsilon^{2}}(\xi_{i})$ depends on $\eta$ too).

Then, by Definition 2.3, we get
$$(KH)\int_{a}^{b}f(t, \omega)d t=C(\omega)\cdot (KH)\int_{a}^{b}\varphi(t) dt,$$
which proves the theorem. $\hfill \square$

{\bf Remark.} Since the Kurzweil-Henstock integrability of a function $\varphi:[a, b]\to \mathbb{R}$ is more general than the Riemann
integrability (in fact, it is equivalent with the so-called Denjoy-Perron integrability, see [2], [8]), Theorem 3.2 gives examples of random functions which are Kurzweil-Henstock integrable in probability on $[a, b]$ but are not Riemann integrable in probability  in the sense of Definition 1.1.

For $p\ge 1$, let us consider
$$L^{p}(E,B,P)=\{g\in L(E,B,P) ; \int_{E}|g(\omega)|^{p}dP(\omega) < +\infty\},$$
where $\int_{E}|g(\omega)|^{q}d P(\omega)$ represents the $q$-th moment of the random variable $g$.

The following Fubini-type result holds.

{\bf Theorem 3.3.} {\it Let $f:[a, b]\to L(E,B,P)$ be Kurzweil-Henstock integrable in probability on $[a, b]$ and such that
there exists $A\in L^{1}(E,B,P)$, $A(\omega)\ge 0$, a.e. $\omega\in E$ with $P(\{\omega\in E; |f(t, \omega)|\le A(\omega)\})=1$, for all
$t\in [a, b]$. Then, $\varphi(t)=\int_{E}f(t, \omega)d P(\omega)$, $t\in [a, b]$, is Kurzweil-Henstock integrable on $[a, b]$ and
$$(KH)\int_{a}^{b}\left [\int_{E}f(t, \omega)d P(\omega)\right ]dt=\int_{E}\left [(KH)\int_{a}^{b}f(t, \omega)d t\right ]d P(\omega).$$}
{\bf Proof.} Let us denote $I(\omega)=(KH)\int_{a}^{b}f(t, \omega)d t\in L(E,B,P)$. Since $f$ is Kurzweil-Henstock integrable on $[a, b]$,
for $\varepsilon >0$ and $\eta=1/m$, $m\in \mathbb{N}$, there exists $\delta_{\varepsilon, m}:[a, b]\to \mathbb{R}$, such that for any division
$d_{m} :a=x_{0}^{(m)} < x_{1}^{(m)} < ... < x_{n_{m}}^{(m)}=b$ and any $\xi_{i}^{(m)}\in [x_{i}^{(m)}, x_{i+1}^{(m)}]$, with
$x_{i+1}^{(m)}-x_{i}^{(m)}<\delta_{\varepsilon, m}(\xi_{i}^{(m)})$, $i=0, ..., n_{m}-1$, we have
$$P(\{\omega\in E; |S(f;d_{m}, \xi_{i}^{(m)})(\omega)-I(\omega)|\ge \varepsilon\})<1/m,\, m\in \mathbb{N}.$$
This means that $S(f;d_{m}, \xi_{i}^{(m)})(\omega)\to I(\omega)$ in probability, as $m\to \infty$.

On the other hand,
$$|S(f;d_{m}, \xi_{i}^{(m)})(\omega)|=|\sum_{i=0}^{n_{m}-1}f(\xi_{i}^{(m)}, \omega)(x^{(m)}_{i+1}-x^{(m)}_{i})|\le A(\omega)\cdot (b-a), \mbox{ a.e. } \omega\in E$$
and taking into account the well-known property of the integral with respect to P, we immediately get
$$\int_{E}I(\omega)dP(\omega)=\lim_{m\to \infty}\int_{E}S(f;d_{m}, \xi_{i}^{(m)})(\omega)d P(\omega)$$
$$=\lim_{m\to \infty}\sum_{i=0}^{n_{m}-1}\left [\int_{E}f(\xi_{i}^{(m)}, \omega)d P(\omega)\right ]\cdot (x_{i+1}^{(m)}-x_{i}^{(m)}).$$
Now, reasoning exactly as in the case of the definitions of Riemann integrability (see e.g. [6, p. 379-380 and p. 383-384]), it is easy to obtain
that the Kurzweil-Henstock integrability in Definition 3.1 is equivalent with the fact that there exists a sequence $\delta_{m}:[a, b]\to \mathbb{R}$, $m\in \mathbb{N}$, such that for any sequence of divisions $(d_{m})_{m\in \mathbb{N}}$, $d_{m}:a=x_{0}^{(m)}<x_{1}^{(m)}< ... <x_{n_{m}}^{(m)}=b$,
and any sequence $(\xi_{i}^{(m)})_{i=0, n_{m}-1}$ with $\xi_{i}^{(m)}\in [x_{i}^{(m)}, x_{i+1}^{(m)}]$, $x_{i+1}^{(m)}-x_{i}^{(m)} < \delta_{m}(\xi_{i}^{(m)})$, $i=0, ..., n_{m}-1$, we have
$$\lim_{m\to \infty}S(\varphi, d_{m}, \xi_{i}^{(m)})=I=(KH)\int_{a}^{b}\varphi(t) dt.$$
But, denoting $\varphi(t)=\int_{E}f(t, \omega)d P(\omega)$, $t\in [a, b]$, by the previous reasonings we immediately get that
$\varphi$ is Kurzweil-Henstock integrable on $[a, b]$ and
$$(KH)\int_{a}^{b}\left [\int_{E}f(t, \omega)d P(\omega)\right ]d t=\int_{E}I(\omega)dP(\omega)$$
$$=\int_{E}\left [(KH)\int_{a}^{b}f(t, \omega)d t\right ]d P(\omega),$$
which proves the theorem. $\hfill \square$

{\bf Remarks.} 1) Theorem 3.3 is an analogue of Theorem III.8 in [7, p. 55].

2) Let $f, F:[a, b]\to \mathbb{R}$ be such that $F^{\prime}(x)=f(x)$, $x\in (a, b)$. It is known (se e.g. [1]) that in this case $f$ is (KH)-integrable on $[a, b]$ and $(KH)\int_{a}^{b}f(x)d x=F(b)-F(a)$.

Now, let $f, F:[a, b]\to L(E,B,P)$ be such that in each $t_{0}\in (a, b)$, $f$ is the derivative in probability of $F(t_{0}, w)$, i.e. for all $\varepsilon, \eta >0$, there exists $\delta(\varepsilon, \eta)>0$, such that for all $t\in [a, b]$, $t\not=t_{0}$, $|t-t_{0}|<\delta(\varepsilon, \eta)$, we have
$$P(\{\omega\in E; | \left [F(t_{0}, \omega)-F(t, \omega)\right ]/(t-t_{0})-f(t_{0}, \omega)|\ge \varepsilon\})<\eta,$$
holds.

Then, the following question arises : in what conditions $f(t, \omega)$ is (KH)-integrable in probability on $[a, b]$ and
$$(KH)\int_{a}^{b}f(t, \omega)d P(\omega)=F(b, \omega)-F(a, \omega), \mbox{ a.e. } \omega\in E.$$

\end{document}